\documentclass[12pt]{article}
\def\Dj{\hbox{D\kern-.73em\raise.30ex\hbox{-}
\raise-.30ex\hbox{}}}
\def\dj{\hbox{d\kern-.33em\raise.80ex\hbox{-}
\raise-.80ex\hbox{\kern-.40em}}}
\usepackage{epsfig}
\usepackage{amsmath,amsthm,amsfonts,amssymb,amscd}
\begin{document}

\baselineskip=0.30in

\parskip=8pt

\newtheorem{lem}{Lemma}[section]
\newtheorem{thm}[lem]{Theorem}
\newtheorem{cor}[lem]{Corollary}
\newtheorem*{prop}{Proposition}
\newtheorem{con}[lem]{Conjecture}
\newtheorem{rem}[lem]{Remark}
\newtheorem{defi}[lem]{Definition}
\renewcommand\baselinestretch{1.1}
\def\pf{\noindent {\it Proof.} }
\def\qed{\hfill \rule{4pt}{7pt}}

\begin{center} {\Large \bf A proof of the conjecture on hypoenergetic\\
graphs with maximum degree $\Delta \leq 3$
 \footnotetext[1]{Supported by NSFC No.10831001, PCSIRT and the ``973" program.}
 \footnotetext[2]{Supported by NSFC No.10871166, NSFJS and NSFUJS.}}
 \end{center}

\begin{center}
{ \small  Xueliang Li\footnotemark[1],~~Hongping Ma\footnotemark[2]\\[5pt]
\small Center for Combinatorics and LPMC-TJKLC,\\
\small Nankai University, Tianjin 300071, P.R. China.\\
\small Email: lxl@nankai.edu.cn; mhp@cfc.nankai.edu.cn}
\end{center}

\begin{abstract}
The energy $E(G)$ of a graph $G$ is defined as the sum of the
absolute values of its eigenvalues. A graph $G$ of order $n$ is said
to be hypoenergetic if $E(G)<n$. Majstorovi\'{c} et al. conjectured
that complete bipartite graph $K_{2,3}$ is the only hypoenergetic
connected quadrangle-containing graph with maximum degree $\Delta
\leq 3$. This paper is devoted to giving a confirmative proof to the
conjecture.\\[2mm]
{\bf Keywords:} energy of a graph; hypoenergetic;
quadrangle-containing (-free); cyclomatic number\\[2mm]
{\bf AMS Subject Classification 2000:} 15A18; 05C50; 05C90; 92E10
\end{abstract}

\section{Introduction}

We use Bondy and Murty \cite{BM} for terminology and notations not
defined here. Let $G$ be a simple graph with $n$ vertices and $m$
edges. The {\it cyclomatic number} of a connected graph $G$ is
defined as $c(G)=m-n+1$. A graph $G$ with $c(G)=k$ is called a {\it
$k$-cyclic graph}. In particular, for $c(G)=0,1,2$ or $3$ we call
$G$ a tree, unicyclic, bicyclic or tricyclic graph, respectively.
Denote by $\Delta$ the maximum degree of a graph. The eigenvalues
$\lambda_{1}, \lambda_{2},\ldots, \lambda_{n}$ of the adjacency
matrix $A(G)$ of $G$ are said to be the eigenvalues of the graph
$G$. The $energy$ of $G$ is defined as
$$E=E(G)=\sum_{i=1}^{n}|\lambda_{i}|.$$

For several classes of graphs it has been demonstrated that the
energy exceeds the number of vertices (see, \cite{G1}). In 2007,
Nikiforov \cite{N} showed that for almost all graphs,
$$E=\left(\frac{4}{3\pi}+o(1)\right)n^{3/2}.$$
Thus the number of graphs satisfying the condition $E<n$ is
relatively small. In \cite{GR}, a $hypoenergetic$ graph is defined
to be a (connected) graph satisfying $E<n$.

Gutman et al. \cite{GLSZ} gave results on hypoenergetic trees. You
and Liu \cite{YL} studied hypoenergetic unicyclic and bicyclic
graphs. You, Liu and Gutman \cite{YLG} considered hypoenergetic
tricyclic and $k$-cyclic graphs. In \cite{LM}, the present authors
showed that there exist hypoenergetic $k$-cyclic graphs of order $n$
and maximum degree $\Delta$ for all (suitable large) $n$ and
$\Delta$; And for $\Delta \geq 4$ there exist hypoenergetic
unicyclic, bicyclic and tricyclic graphs for all $n$ except very few
small values of $n$. For hypoenergetic graphs with $\Delta \leq 3$,
we have the following results.

\begin{figure}[ht]
\centering
  \setlength{\unitlength}{0.05 mm}%
  \begin{picture}(1592.3, 330.9)(0,0)
  \put(0,0){\includegraphics{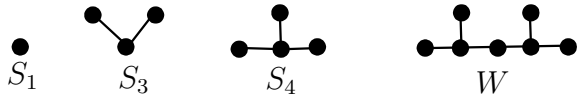}}
  \put(40.00,60.93){\fontsize{11.76}{14.11}\selectfont \makebox(206.6, 82.6)[l]{$S_1$\strut}}
  \put(335.53,52.15){\fontsize{11.76}{14.11}\selectfont \makebox(206.6, 82.6)[l]{$S_3$\strut}}
  \put(727.62,49.22){\fontsize{11.76}{14.11}\selectfont \makebox(206.6, 82.6)[l]{$S_4$\strut}}
  \put(1286.50,37.52){\fontsize{11.76}{14.11}\selectfont \makebox(124.0, 82.6)[l]{$W$\strut}}
  \end{picture}%
 \caption{The hypoenergetic trees with maximum degree at most $3$.} \label{fig1}
\end{figure}

\begin{lem}\cite{GLSZ}\label{lem1.1}
There exist only four hypoenergetic trees with $\Delta \leq 3$,
dipicted in Figure \ref{fig1}.
\end{lem}

\begin{lem}\cite{N2}\label{lem1.2}
Let $G$ be a graph of order $n$ with at least $n$ edges and with no
isolated vertices. If $G$ is quadrangle-free and $\Delta (G)\leq 3$,
then $E(G)>n$.
\end{lem}

In \cite{MKG} Majstorovi\'{c} et al. proposed the following
conjecture, which is the first half of their Conjecture 3.7.

\begin{con}\cite{MKG}\label{con1.3}
Complete bipartite graph $K_{2,3}$ is the only hypoenergetic
connected quadrangle-containing graph with $\Delta \leq 3$.
\end{con}

It follows from Lemma \ref{lem1.2} that Conjecture \ref{con1.3} is
equivalent to the following result.

\begin{thm}\label{thm1.4}
$K_{2,3}$ is the only hypoenergetic connected cyclic graph with
$\Delta \leq 3$.
\end{thm}

We will give a proof of Theorem \ref{thm1.4} in the next section.
Therefore, combining Lemma \ref{lem1.1}, we obtain

\begin{thm}\label{thm1.5}
$S_1, S_3, S_4, W$ (see Figure \ref{fig1}) and $K_{2,3}$ are the
only hypoenergetic connected graphs with $\Delta \leq 3$.
\end{thm}

\section{Main results}

The following two lemmas are need in the sequel.

\begin{lem}\cite{LM}\label{lem2.1}
$K_{2,3}$ is the only hypoenergetic graph with $\Delta \leq 3$ among
all unicyclic and bicyclic graphs.
\end{lem}

\begin{lem}\cite{DS}\label{lem2.2}
If $F$ is an edge cut of a simple graph $G$, then $E(G-F)\leq E(G)$,
where $G-F$ is the subgraph obtained from $G$ by deleting the edges
in $F$.
\end{lem}

\noindent{\bf Proof of Theorem \ref{thm1.4}:} Notice that $K_{2,3}$
is hypoenergetic by Lemma \ref{lem2.1}. Let $G$ be a connected
cyclic graph with $G\not\cong K_{2,3}$, $\Delta \leq 3$ and
$c(G)=m-n+1\geq 1$. In the following we show that $G$ is
non-hypoenergetic by induction on $c(G)$. It follows from Lemma
\ref{lem2.1} that the result is true if $c(G)\leq 2$. We assume that
$G$ is non-hypoenergetic for $1\leq c(G)< k$. Now let $G$ be a graph
with $c(G)=k\geq 3$. In the following we will repeatedly make use of
the following claim:

{\bf Claim 1.} {\it If there exists an edge cut $F$ of $G$ such that
$G-F$ has exactly two components $G_1$, $G_2$ with $0\leq c(G_1),
c(G_2)<k$ and $G_1, G_2\not\cong S_1, S_3, S_4, W, K_{2,3}$, then we
are done.}

\pf It follows from Lemma \ref{lem1.1} and the induction hypothesis
that $G_1$ and $G_2$ are non-hypoenergetic. By Lemma \ref{lem2.2},
we have $E(G)\geq E(G-F)$. Therefore
$$E(G)\geq E(G-F)=E(G_1)+E(G_2)\geq |V(G_1)|+|V(G_2)|=n,$$
which proves the claim. \qed

For convenience, we call an edge cut $F$ of $G$ a {\it good edge
cut} if $F$ satisfies the conditions in Claim 1. In what follows, we
use $\bar{G}$ to denote the graph obtained from $G$ by repeatedly
deleting the pendent vertices. Clearly, $c(\bar{G})=c(G)$. Denote by
$\kappa'(\bar{G})$ the edge connectivity of $\bar{G}$. Since
$\Delta(\bar{G})\leq 3$, we have $1\leq \kappa'(\bar{G})\leq 3$.
Therefore, we only need to consider the following three cases.

{\bf Case 1.} $\kappa'(\bar{G})=1$.

\begin{figure}[ht]
\centering
    \setlength{\unitlength}{0.05 mm}%
  \begin{picture}(1942.8, 556.4)(0,0)
  \put(0,0){\includegraphics{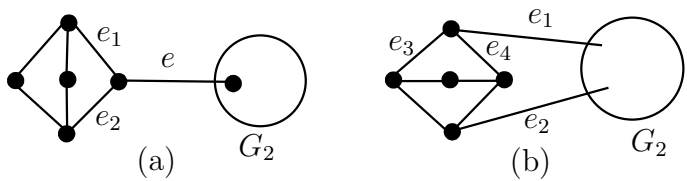}}
  \put(446.11,312.33){\fontsize{11.76}{14.11}\selectfont \makebox(124.0, 82.6)[l]{$e$\strut}}
  \put(275.16,383.07){\fontsize{11.76}{14.11}\selectfont \makebox(206.6, 82.6)[l]{$e_1$\strut}}
  \put(272.21,159.07){\fontsize{11.76}{14.11}\selectfont \makebox(206.6, 82.6)[l]{$e_2$\strut}}
  \put(653.30,77.44){\fontsize{11.76}{14.11}\selectfont \makebox(206.6, 82.6)[l]{$G_2$\strut}}
  \put(1051.57,357.38){\fontsize{11.76}{14.11}\selectfont \makebox(206.6, 82.6)[l]{$e_3$\strut}}
  \put(1421.83,431.27){\fontsize{11.76}{14.11}\selectfont \makebox(206.6, 82.6)[l]{$e_1$\strut}}
  \put(1411.43,150.39){\fontsize{11.76}{14.11}\selectfont \makebox(206.6, 82.6)[l]{$e_2$\strut}}
  \put(1696.20,89.92){\fontsize{11.76}{14.11}\selectfont \makebox(206.6, 82.6)[l]{$G_2$\strut}}
  \put(1309.64,345.79){\fontsize{11.76}{14.11}\selectfont \makebox(206.6, 82.6)[l]{$e_4$\strut}}
  \put(381.35,40.02){\fontsize{11.76}{14.11}\selectfont \makebox(124.0, 82.6)[l]{(a)\strut}}
  \put(1374.35,37.52){\fontsize{11.76}{14.11}\selectfont \makebox(124.0, 82.6)[l]{(b)\strut}}
  \end{picture}%
  \caption{The graphs in Case 1 and Subcase 2.1 of Theorem \ref{thm1.4}.} \label{fig2}
\end{figure}

Let $e$ be a cut edge of $\bar{G}$. Then $\bar{G}-e$ has exactly two
components, say, $H_1$ and $H_2$. It is clear that $c(H_1)\geq 1$,
$c(H_2)\geq 1$ and $c(H_1)+c(H_2)=k$. Consequently, $G-e$ has
exactly two components $G_1$ and $G_2$ with $c(G_1)\geq 1$,
$c(G_2)\geq 1$ and $c(G_1)+c(G_2)=k$, where $H_i$ is a subgraph of
$G_i$ for $i=1,2$. If neither $G_1$ nor $G_2$ is isomorphic to
$K_{2,3}$, then we are done by Claim 1. Otherwise, by symmetry we
assume that $G_1\cong K_{2,3}$. Then $G$ must have the structure as
given in Figure \ref{fig2} (a). Now, let $F=\{e_1,e_2\}$. Then $G-F$
has exactly two components $G'_1$ and $G'_2$, where $G'_1$ is a
quadrangle and $G'_2$ is a graph obtained from $G_2$ by adding a
pendent edge. Therefore we have that $c(G'_2)=k-2$ and
$G'_2\not\cong K_{2,3}$, and so we are done by Claim 1.

{\bf Case 2.} $\kappa'(\bar{G})=2$.

Let $F=\{e_1, e_2\}$ be an edge cut of $\bar{G}$. Then $\bar{G}-F$
has exactly two components, say, $H_1$ and $H_2$. Clearly,
$c(H_1)+c(H_2)=k-1\geq 2$.

{\bf Subcase 2.1.} $c(H_1)\geq 1$ and $c(H_2)\geq 1$.

Therefore, $G-F$ has exactly two components $G_1$ and $G_2$ with
$c(G_1)\geq 1$, $c(G_2)\geq 1$ and $c(G_1)+c(G_2)=k-1$, where $H_i$
is a subgraph of $G_i$ for $i=1,2$. If neither $G_1$ nor $G_2$ is
isomorphic to $K_{2,3}$, then we are done by Claim 1. Otherwise, by
symmetry we assume that $G_1\cong K_{2,3}$. Then $G$ must have the
structure as given in Figure \ref{fig2} (b). Now, let
$F'=\{e_2,e_3,e_4\}$. Then it is easy to see that $F'$ is a good
edge cut. The proof is thus complete.

{\bf Subcase 2.2.} One of $H_1$ and $H_2$, say $H_2$ is a tree.

Therefore, $G-F$ has exactly two components $G_1$ and $G_2$ with
$c(G_1)=k-1$ and $c(G_2)=0$, where $H_i$ is a subgraph of $G_i$ for
$i=1,2$. If $G_1\not\cong K_{2,3}$ and $G_2\not\cong S_1,S_3,S_4,W$,
then we are done by Claim 1. So we assume that this is not true. We
only need to consider the following five cases.

{\bf Subsubcase 2.2.1.} $G_2\cong S_1$.

\begin{figure}[ht]
\centering
    \setlength{\unitlength}{0.05 mm}%
  \begin{picture}(2102.3, 1323.0)(0,0)
  \put(0,0){\includegraphics{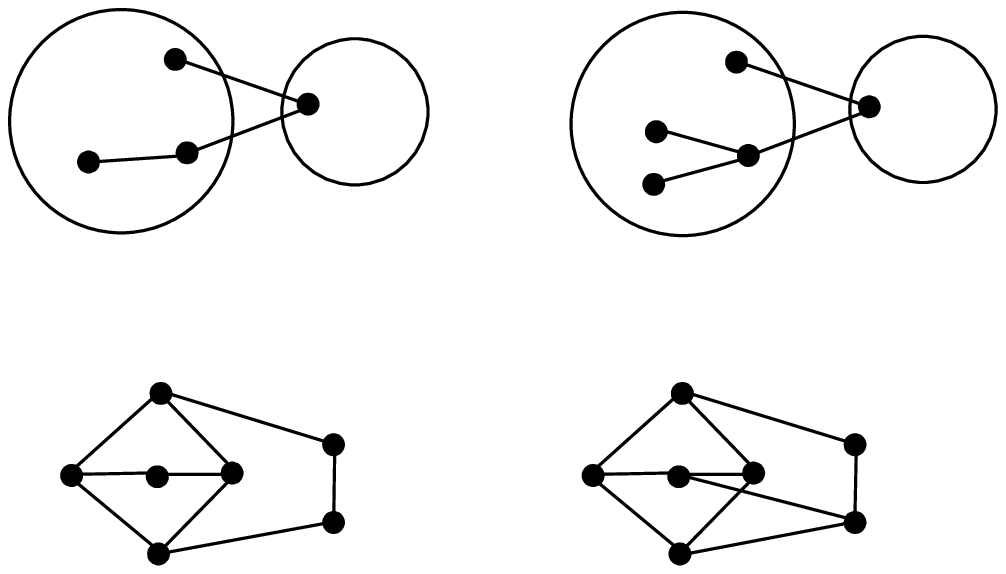}}
  \put(681.46,1063.80){\fontsize{11.76}{14.11}\selectfont \makebox(124.0, 82.6)[l]{$x$\strut}}
  \put(268.81,1155.45){\fontsize{11.76}{14.11}\selectfont \makebox(206.6, 82.6)[l]{$x_1$\strut}}
  \put(368.43,900.52){\fontsize{11.76}{14.11}\selectfont \makebox(206.6, 82.6)[l]{$x_2$\strut}}
  \put(512.66,1124.26){\fontsize{11.76}{14.11}\selectfont \makebox(206.6, 82.6)[l]{$e_1$\strut}}
  \put(517.70,975.62){\fontsize{11.76}{14.11}\selectfont \makebox(206.6, 82.6)[l]{$e_2$\strut}}
  \put(1830.28,1063.80){\fontsize{11.76}{14.11}\selectfont \makebox(124.0, 82.6)[l]{$x$\strut}}
  \put(1414.34,1167.72){\fontsize{11.76}{14.11}\selectfont \makebox(206.6, 82.6)[l]{$x_1$\strut}}
  \put(1508.99,895.18){\fontsize{11.76}{14.11}\selectfont \makebox(206.6, 82.6)[l]{$x_2$\strut}}
  \put(1658.97,1131.82){\fontsize{11.76}{14.11}\selectfont \makebox(206.6, 82.6)[l]{$e_1$\strut}}
  \put(1656.45,973.10){\fontsize{11.76}{14.11}\selectfont \makebox(206.6, 82.6)[l]{$e_2$\strut}}
  \put(424.31,622.73){\fontsize{11.76}{14.11}\selectfont \makebox(124.0, 82.6)[l]{(a)\strut}}
  \put(1620.96,625.40){\fontsize{11.76}{14.11}\selectfont \makebox(124.0, 82.6)[l]{(b)\strut}}
  \put(245.61,731.24){\fontsize{11.76}{14.11}\selectfont \makebox(206.6, 82.6)[l]{$G_1$\strut}}
  \put(1391.92,723.68){\fontsize{11.76}{14.11}\selectfont \makebox(206.6, 82.6)[l]{$G_1$\strut}}
  \put(685.29,816.21){\fontsize{11.76}{14.11}\selectfont \makebox(206.6, 82.6)[l]{$G_2$\strut}}
  \put(1855.64,823.08){\fontsize{11.76}{14.11}\selectfont \makebox(206.6, 82.6)[l]{$G_2$\strut}}
  \put(331.71,529.27){\fontsize{11.76}{14.11}\selectfont \makebox(206.6, 82.6)[l]{$x_1$\strut}}
  \put(729.97,199.66){\fontsize{11.76}{14.11}\selectfont \makebox(206.6, 82.6)[l]{$x_2$\strut}}
  \put(175.65,887.55){\fontsize{11.76}{14.11}\selectfont \makebox(206.6, 82.6)[l]{$y_1$\strut}}
  \put(1559.07,37.52){\fontsize{11.76}{14.11}\selectfont \makebox(124.0, 82.6)[l]{(d)\strut}}
  \put(1341.53,837.06){\fontsize{11.76}{14.11}\selectfont \makebox(206.6, 82.6)[l]{$y_1$\strut}}
  \put(1326.41,1056.24){\fontsize{11.76}{14.11}\selectfont \makebox(206.6, 82.6)[l]{$y_2$\strut}}
  \put(480.77,454.60){\fontsize{11.76}{14.11}\selectfont \makebox(206.6, 82.6)[l]{$e_1$\strut}}
  \put(721.07,400.70){\fontsize{11.76}{14.11}\selectfont \makebox(124.0, 82.6)[l]{$x$\strut}}
  \put(720.15,282.57){\fontsize{11.76}{14.11}\selectfont \makebox(206.6, 82.6)[l]{$e_2$\strut}}
  \put(425.80,39.97){\fontsize{11.76}{14.11}\selectfont \makebox(124.0, 82.6)[l]{(c)\strut}}
  \put(352.21,388.05){\fontsize{11.76}{14.11}\selectfont \makebox(206.6, 82.6)[l]{$e_3$\strut}}
  \put(241.83,334.08){\fontsize{11.76}{14.11}\selectfont \makebox(206.6, 82.6)[l]{$e_4$\strut}}
  \put(185.41,204.08){\fontsize{11.76}{14.11}\selectfont \makebox(206.6, 82.6)[l]{$e_5$\strut}}
  \put(1391.39,529.27){\fontsize{11.76}{14.11}\selectfont \makebox(206.6, 82.6)[l]{$x_1$\strut}}
  \put(1789.64,199.66){\fontsize{11.76}{14.11}\selectfont \makebox(206.6, 82.6)[l]{$x_2$\strut}}
  \put(1540.45,454.60){\fontsize{11.76}{14.11}\selectfont \makebox(206.6, 82.6)[l]{$e_1$\strut}}
  \put(1780.75,400.70){\fontsize{11.76}{14.11}\selectfont \makebox(124.0, 82.6)[l]{$x$\strut}}
  \put(1779.83,282.57){\fontsize{11.76}{14.11}\selectfont \makebox(206.6, 82.6)[l]{$e_2$\strut}}
  \put(1411.89,388.05){\fontsize{11.76}{14.11}\selectfont \makebox(206.6, 82.6)[l]{$e_3$\strut}}
  \put(1301.50,334.08){\fontsize{11.76}{14.11}\selectfont \makebox(206.6, 82.6)[l]{$e_4$\strut}}
  \put(1245.09,204.08){\fontsize{11.76}{14.11}\selectfont \makebox(206.6, 82.6)[l]{$e_5$\strut}}
  \end{picture}%
  \caption{The graphs in Subsubcase 2.2.1 of Theorem \ref{thm1.4}.} \label{fig3}
\end{figure}

Let $V(G_2)=\{x\}$, $e_1=xx_1$ and $e_2=xx_2$. It is clear that
$d_{G_1}(x_2)=1$ or $2$. If $d_{G_1}(x_2)=1$, let
$N_{G_1}(x_2)=\{y_1\}$ (see Figure \ref{fig3} (a), where $y_1$ may
be equal to $x_1$). Let $F'=\{e_1, x_2y_1\}$. Then $G-F'$ has
exactly two components $G'_1$ and $G'_2$, where $G'_1$ is a graph
obtained from $G_1$ by deleting a pendent vertex and $G'_2$ is a
tree of order $2$. Therefore, $c(G'_1)=k-1$. If $G'_1\not\cong
K_{2,3}$, then we are done by Claim 1. Otherwise, $G$ must be the
graph as given in Figure \ref{fig3} (c). It is easy to see that
$F''=\{e_1,e_3,e_4,e_5\}$ is a good edge cut.

If $d_{G_1}(x_2)=2$, let $N_{G_1}(x_2)=\{y_1,y_2\}$ (see Figure
\ref{fig3} (b), where one of $y_1$ and $y_2$ may be equal to $x_1$).
Let $F'=\{e_1, x_2y_1, x_2y_2\}$. Then $G-F'$ has exactly two
components $G'_1$ and $G'_2$ such that $G'_1$ is a graph obtained
from $G_1$ by deleting a vertex of degree $2$ and $G'_2$ is a tree
of order $2$. Therefore, $c(G'_1)=k-2$. If $G'_1\not\cong K_{2,3}$,
then we are done by Claim 1. Otherwise, $G$ must be the graph as
given in Figure \ref{fig3} (d). It is easy to see that
$F''=\{e_1,e_3,e_4,e_5\}$ is a good edge cut.

{\bf Subsubcase 2.2.2.} $G_2\cong S_3$.

\begin{figure}[ht]
\centering
 \setlength{\unitlength}{0.05 mm}%
  \begin{picture}(2492.5, 1244.6)(0,0)
  \put(0,0){\includegraphics{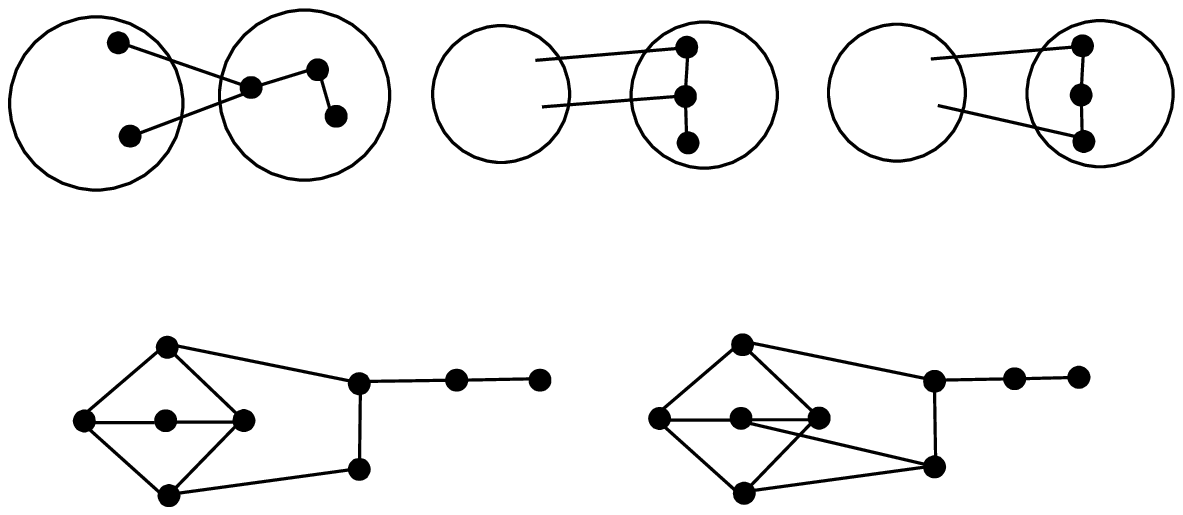}}
  \put(528.29,1076.95){\fontsize{11.76}{14.11}\selectfont \makebox(124.0, 82.6)[l]{$x$\strut}}
  \put(153.15,1087.67){\fontsize{11.76}{14.11}\selectfont \makebox(206.6, 82.6)[l]{$x_1$\strut}}
  \put(145.11,902.78){\fontsize{11.76}{14.11}\selectfont \makebox(206.6, 82.6)[l]{$x_2$\strut}}
  \put(396.92,1076.80){\fontsize{11.76}{14.11}\selectfont \makebox(206.6, 82.6)[l]{$e_1$\strut}}
  \put(401.96,928.16){\fontsize{11.76}{14.11}\selectfont \makebox(206.6, 82.6)[l]{$e_2$\strut}}
  \put(386.64,650.50){\fontsize{11.76}{14.11}\selectfont \makebox(124.0, 82.6)[l]{(a)\strut}}
  \put(206.74,739.32){\fontsize{11.76}{14.11}\selectfont \makebox(206.6, 82.6)[l]{$G_1$\strut}}
  \put(579.20,752.72){\fontsize{11.76}{14.11}\selectfont \makebox(206.6, 82.6)[l]{$G_2$\strut}}
  \put(348.76,534.54){\fontsize{12.93}{15.52}\selectfont \makebox(227.3, 90.9)[l]{$x_1$\strut}}
  \put(725.77,468.90){\fontsize{12.93}{15.52}\selectfont \makebox(136.4, 90.9)[l]{$x$\strut}}
  \put(768.03,219.09){\fontsize{12.93}{15.52}\selectfont \makebox(227.3, 90.9)[l]{$x_2$\strut}}
  \put(1506.33,545.26){\fontsize{12.93}{15.52}\selectfont \makebox(227.3, 90.9)[l]{$x_1$\strut}}
  \put(1886.28,471.87){\fontsize{12.93}{15.52}\selectfont \makebox(136.4, 90.9)[l]{$x$\strut}}
  \put(1937.01,224.13){\fontsize{12.93}{15.52}\selectfont \makebox(227.3, 90.9)[l]{$x_2$\strut}}
  \put(1642.68,42.92){\fontsize{11.76}{14.11}\selectfont \makebox(124.0, 82.6)[l]{(e)\strut}}
  \put(514.65,37.52){\fontsize{11.76}{14.11}\selectfont \makebox(124.0, 82.6)[l]{(d)\strut}}
  \put(1210.58,1116.59){\fontsize{11.76}{14.11}\selectfont \makebox(206.6, 82.6)[l]{$e_1$\strut}}
  \put(1213.28,943.90){\fontsize{11.76}{14.11}\selectfont \makebox(206.6, 82.6)[l]{$e_2$\strut}}
  \put(1181.12,655.19){\fontsize{11.76}{14.11}\selectfont \makebox(124.0, 82.6)[l]{(b)\strut}}
  \put(1026.69,792.92){\fontsize{11.76}{14.11}\selectfont \makebox(206.6, 82.6)[l]{$G_1$\strut}}
  \put(1394.28,790.10){\fontsize{11.76}{14.11}\selectfont \makebox(206.6, 82.6)[l]{$G_2$\strut}}
  \put(2014.72,1119.45){\fontsize{11.76}{14.11}\selectfont \makebox(206.6, 82.6)[l]{$e_1$\strut}}
  \put(2001.23,906.29){\fontsize{11.76}{14.11}\selectfont \makebox(206.6, 82.6)[l]{$e_2$\strut}}
  \put(1985.25,658.05){\fontsize{11.76}{14.11}\selectfont \makebox(124.0, 82.6)[l]{(c)\strut}}
  \put(1826.06,795.66){\fontsize{11.76}{14.11}\selectfont \makebox(206.6, 82.6)[l]{$G_1$\strut}}
  \put(2197.98,787.57){\fontsize{11.76}{14.11}\selectfont \makebox(206.6, 82.6)[l]{$G_2$\strut}}
  \put(1435.67,1046.92){\fontsize{11.76}{14.11}\selectfont \makebox(206.6, 82.6)[l]{$e_3$\strut}}
  \put(2245.89,1044.79){\fontsize{11.76}{14.11}\selectfont \makebox(206.6, 82.6)[l]{$e_3$\strut}}
  \put(525.61,470.23){\fontsize{12.93}{15.52}\selectfont \makebox(227.3, 90.9)[l]{$e_1$\strut}}
  \put(1688.54,486.31){\fontsize{12.93}{15.52}\selectfont \makebox(227.3, 90.9)[l]{$e_1$\strut}}
  \put(196.02,231.75){\fontsize{12.93}{15.52}\selectfont \makebox(227.3, 90.9)[l]{$e_3$\strut}}
  \put(1380.39,237.11){\fontsize{12.93}{15.52}\selectfont \makebox(227.3, 90.9)[l]{$e_3$\strut}}
  \put(429.27,238.58){\fontsize{12.93}{15.52}\selectfont \makebox(227.3, 90.9)[l]{$e_4$\strut}}
  \put(1602.92,247.69){\fontsize{12.93}{15.52}\selectfont \makebox(227.3, 90.9)[l]{$e_4$\strut}}
  \put(1742.13,296.06){\fontsize{12.93}{15.52}\selectfont \makebox(227.3, 90.9)[l]{$e_5$\strut}}
  \end{picture}%
 \caption{The graphs in Subsubcase 2.2.2 of Theorem \ref{thm1.4}.} \label{fig4}
\end{figure}

If $e_1$, $e_2$ are incident with a common vertex in $G_2$, then $G$
must have the structure as given in Figure \ref{fig4} (a). Similar
to the proof of Subsubcase 2.2.1, we can obtain that there exists an
edge cut $F'$ such that $G-F'$ has exactly two components $G'_1$ and
$G'_2$ satisfying that $c(G'_1)=k-1$ if $d_{G_1}(x_2)=1$ or
$c(G'_1)=k-2$ if $d_{G_1}(x_2)=2$ and $G'_2$ is a path of order $4$.
If $G'_1\not\cong K_{2,3}$, then we are done by Claim 1. Otherwise
$G$ must be the graph as given in Figure \ref{fig4} (d) or (e). In
the former case $F''=\{e_1,e_3,e_4\}$ is a good edge cut while in
the latter case $F''=\{e_1,e_3,e_4,e_5\}$ is a good edge cut.

If $e_1$, $e_2$ are incident with two different vertices in $G_2$,
then $G$ must have the structure as given in Figure \ref{fig4} (b)
or (c). It is easy to see that $F'=\{e_2, e_3\}$ is a good edge cut.
The proof is thus complete.

{\bf Subsubcase 2.2.3.} $G_2\cong S_4$.

\begin{figure}[ht]
\centering
  \setlength{\unitlength}{0.05 mm}%
  \begin{picture}(2064.3, 1168.4)(0,0)
  \put(0,0){\includegraphics{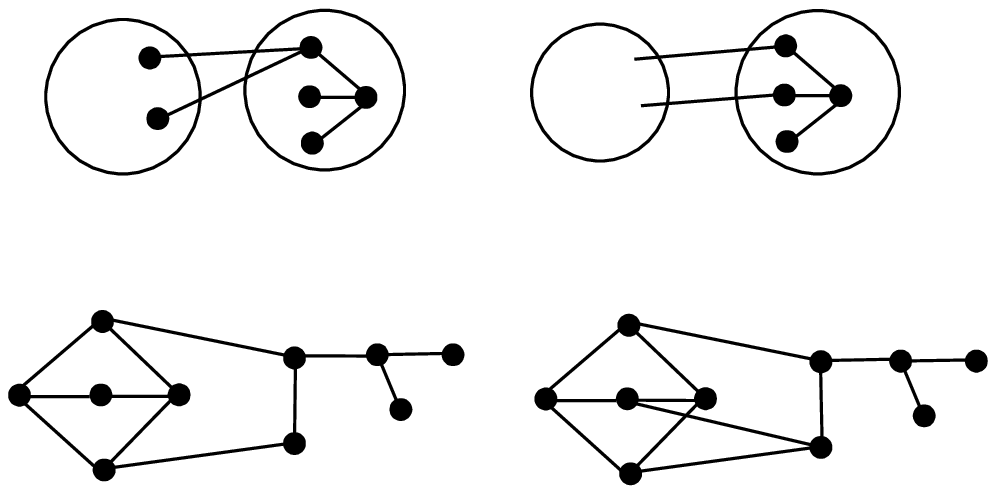}}
  \put(215.03,981.09){\fontsize{11.76}{14.11}\selectfont \makebox(206.6, 82.6)[l]{$x_1$\strut}}
  \put(215.03,863.19){\fontsize{11.76}{14.11}\selectfont \makebox(206.6, 82.6)[l]{$x_2$\strut}}
  \put(444.00,575.26){\fontsize{11.76}{14.11}\selectfont \makebox(124.0, 82.6)[l]{(a)\strut}}
  \put(681.27,1029.32){\fontsize{11.76}{14.11}\selectfont \makebox(124.0, 82.6)[l]{$x$\strut}}
  \put(1483.73,37.52){\fontsize{11.76}{14.11}\selectfont \makebox(124.0, 82.6)[l]{(d)\strut}}
  \put(423.72,37.52){\fontsize{11.76}{14.11}\selectfont \makebox(124.0, 82.6)[l]{(c)\strut}}
  \put(1411.58,1040.08){\fontsize{11.76}{14.11}\selectfont \makebox(206.6, 82.6)[l]{$e_1$\strut}}
  \put(1414.28,867.39){\fontsize{11.76}{14.11}\selectfont \makebox(206.6, 82.6)[l]{$e_2$\strut}}
  \put(1382.12,578.67){\fontsize{11.76}{14.11}\selectfont \makebox(124.0, 82.6)[l]{(b)\strut}}
  \put(1227.90,700.32){\fontsize{11.76}{14.11}\selectfont \makebox(206.6, 82.6)[l]{$G_1$\strut}}
  \put(1597.68,692.28){\fontsize{11.76}{14.11}\selectfont \makebox(206.6, 82.6)[l]{$G_2$\strut}}
  \put(446.94,1036.81){\fontsize{11.76}{14.11}\selectfont \makebox(206.6, 82.6)[l]{$e_1$\strut}}
  \put(449.64,864.12){\fontsize{11.76}{14.11}\selectfont \makebox(206.6, 82.6)[l]{$e_2$\strut}}
  \put(260.58,689.02){\fontsize{11.76}{14.11}\selectfont \makebox(206.6, 82.6)[l]{$G_1$\strut}}
  \put(635.72,699.74){\fontsize{11.76}{14.11}\selectfont \makebox(206.6, 82.6)[l]{$G_2$\strut}}
  \put(1646.13,1043.31){\fontsize{11.76}{14.11}\selectfont \makebox(124.0, 82.6)[l]{$x$\strut}}
  \put(1771.85,928.09){\fontsize{11.76}{14.11}\selectfont \makebox(124.0, 82.6)[l]{$y$\strut}}
  \put(1584.29,957.56){\fontsize{11.76}{14.11}\selectfont \makebox(124.0, 82.6)[l]{$z$\strut}}
  \put(594.07,442.08){\fontsize{12.93}{15.52}\selectfont \makebox(136.4, 90.9)[l]{$x$\strut}}
  \put(636.34,192.28){\fontsize{12.93}{15.52}\selectfont \makebox(227.3, 90.9)[l]{$x_2$\strut}}
  \put(1697.20,441.03){\fontsize{12.93}{15.52}\selectfont \makebox(136.4, 90.9)[l]{$x$\strut}}
  \put(1706.02,184.50){\fontsize{12.93}{15.52}\selectfont \makebox(227.3, 90.9)[l]{$x_2$\strut}}
  \put(393.92,443.41){\fontsize{12.93}{15.52}\selectfont \makebox(227.3, 90.9)[l]{$e_1$\strut}}
  \put(1472.98,437.83){\fontsize{12.93}{15.52}\selectfont \makebox(227.3, 90.9)[l]{$e_1$\strut}}
  \put(64.33,204.93){\fontsize{12.93}{15.52}\selectfont \makebox(227.3, 90.9)[l]{$e_3$\strut}}
  \put(1149.40,197.48){\fontsize{12.93}{15.52}\selectfont \makebox(227.3, 90.9)[l]{$e_3$\strut}}
  \put(297.57,211.76){\fontsize{12.93}{15.52}\selectfont \makebox(227.3, 90.9)[l]{$e_4$\strut}}
  \put(1371.93,208.06){\fontsize{12.93}{15.52}\selectfont \makebox(227.3, 90.9)[l]{$e_4$\strut}}
  \put(1511.14,256.43){\fontsize{12.93}{15.52}\selectfont \makebox(227.3, 90.9)[l]{$e_5$\strut}}
  \end{picture}%
  \caption{The graphs in Subsubcase 2.2.3 of Theorem \ref{thm1.4}.} \label{fig5}
\end{figure}

If $e_1$, $e_2$ are incident with a common vertex in $G_2$, then $G$
must have the structure as given in Figure \ref{fig5} (a). Similar
to the proof of Subsubcase 2.2.1, we can obtain that
there exists an edge cut $F'$ such that $G-F'$ has exactly two
components $G'_1$ and $G'_2$ satisfying that $c(G'_1)=k-1$ if
$d_{G_1}(x_2)=1$ or $c(G'_1)=k-2$ if $d_{G_1}(x_2)=2$ and $G'_2$ is
a tree of order $5$. If $G'_1\not\cong K_{2,3}$, then we are done by
Claim 1. Otherwise $G$ is the graph as given in Figure \ref{fig5}
(c) or (d). In the former case $F''=\{e_1,e_3,e_4\}$ is a good edge
cut while in the latter case $F''=\{e_1,e_3,e_4,e_5\}$ is a good
edge cut.

If $e_1$, $e_2$ are incident with two different vertices in $G_2$,
then $G$ must have the structure as given in Figure \ref{fig5} (b).
It is easy to see that $F'=\{xy, yz\}$ is a good edge cut. The proof
is thus complete.

{\bf Subsubcase 2.2.4.} $G_2\cong W$.

\begin{figure}[ht]
\centering
  \setlength{\unitlength}{0.05 mm}%
  \begin{picture}(2289.2, 2088.3)(0,0)
  \put(0,0){\includegraphics{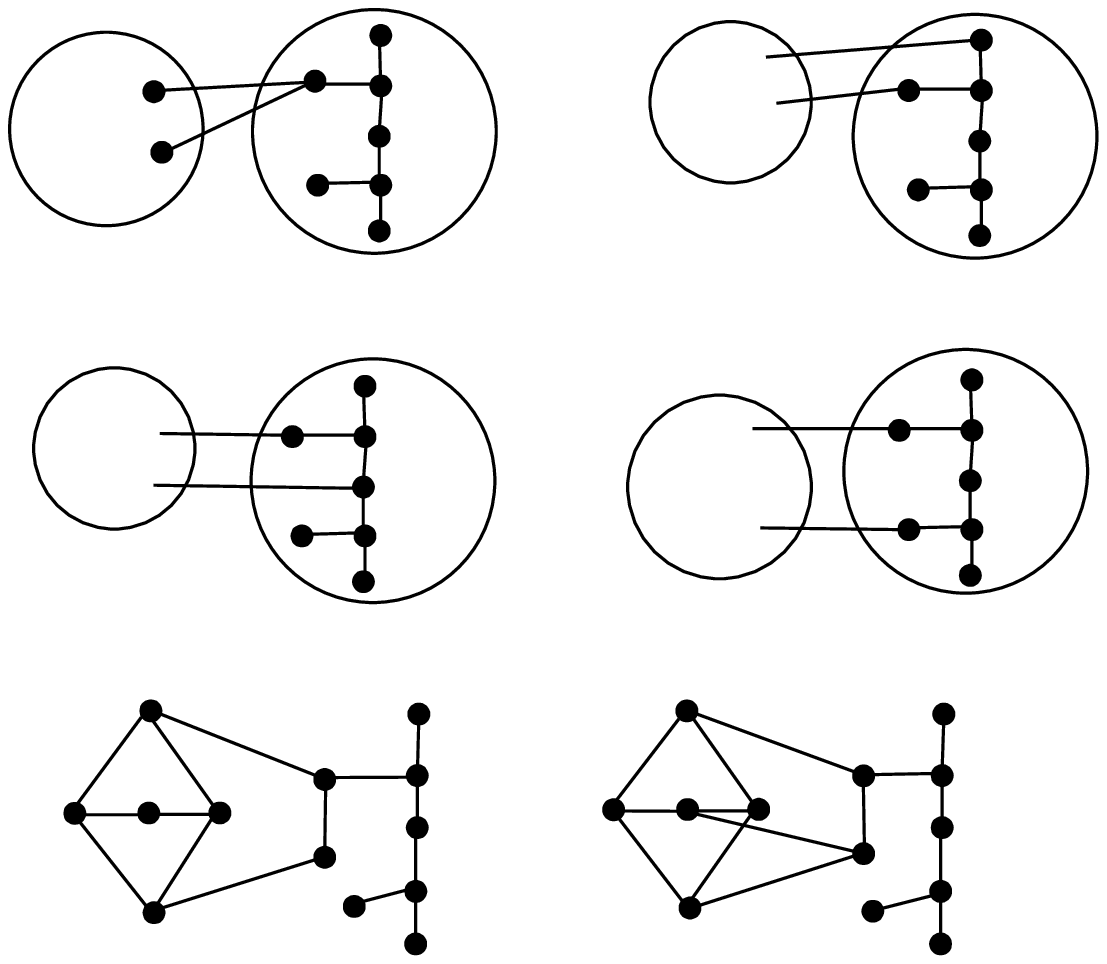}}
  \put(223.32,1829.09){\fontsize{11.76}{14.11}\selectfont \makebox(206.6, 82.6)[l]{$x_1$\strut}}
  \put(223.32,1711.18){\fontsize{11.76}{14.11}\selectfont \makebox(206.6, 82.6)[l]{$x_2$\strut}}
  \put(457.72,1349.30){\fontsize{11.76}{14.11}\selectfont \makebox(124.0, 82.6)[l]{(a)\strut}}
  \put(689.56,1877.32){\fontsize{11.76}{14.11}\selectfont \makebox(124.0, 82.6)[l]{$x$\strut}}
  \put(1768.32,511.08){\fontsize{11.76}{14.11}\selectfont \makebox(124.0, 82.6)[l]{$x$\strut}}
  \put(1587.59,52.67){\fontsize{11.76}{14.11}\selectfont \makebox(124.0, 82.6)[l]{(f)\strut}}
  \put(479.11,37.52){\fontsize{11.76}{14.11}\selectfont \makebox(124.0, 82.6)[l]{(e)\strut}}
  \put(1683.14,1366.11){\fontsize{11.76}{14.11}\selectfont \makebox(124.0, 82.6)[l]{(b)\strut}}
  \put(460.67,1893.16){\fontsize{11.76}{14.11}\selectfont \makebox(206.6, 82.6)[l]{$e_1$\strut}}
  \put(473.48,1742.62){\fontsize{11.76}{14.11}\selectfont \makebox(206.6, 82.6)[l]{$e_2$\strut}}
  \put(226.84,1489.57){\fontsize{11.76}{14.11}\selectfont \makebox(206.6, 82.6)[l]{$G_1$\strut}}
  \put(716.92,1425.51){\fontsize{11.76}{14.11}\selectfont \makebox(206.6, 82.6)[l]{$G_2$\strut}}
  \put(667.27,501.94){\fontsize{11.76}{14.11}\selectfont \makebox(124.0, 82.6)[l]{$x$\strut}}
  \put(1719.49,1963.14){\fontsize{11.76}{14.11}\selectfont \makebox(206.6, 82.6)[l]{$e_1$\strut}}
  \put(1684.26,1780.25){\fontsize{11.76}{14.11}\selectfont \makebox(206.6, 82.6)[l]{$e_2$\strut}}
  \put(1479.26,1575.89){\fontsize{11.76}{14.11}\selectfont \makebox(206.6, 82.6)[l]{$G_1$\strut}}
  \put(1937.30,1425.51){\fontsize{11.76}{14.11}\selectfont \makebox(206.6, 82.6)[l]{$G_2$\strut}}
  \put(1855.56,1893.62){\fontsize{11.76}{14.11}\selectfont \makebox(124.0, 82.6)[l]{$x$\strut}}
  \put(694.50,719.57){\fontsize{11.76}{14.11}\selectfont \makebox(206.6, 82.6)[l]{$G_2$\strut}}
  \put(413.09,697.64){\fontsize{11.76}{14.11}\selectfont \makebox(124.0, 82.6)[l]{(c)\strut}}
  \put(463.87,1192.67){\fontsize{11.76}{14.11}\selectfont \makebox(206.6, 82.6)[l]{$e_1$\strut}}
  \put(447.86,988.47){\fontsize{11.76}{14.11}\selectfont \makebox(206.6, 82.6)[l]{$e_2$\strut}}
  \put(220.44,860.50){\fontsize{11.76}{14.11}\selectfont \makebox(206.6, 82.6)[l]{$G_1$\strut}}
  \put(1927.70,741.99){\fontsize{11.76}{14.11}\selectfont \makebox(206.6, 82.6)[l]{$G_2$\strut}}
  \put(1699.06,677.74){\fontsize{11.76}{14.11}\selectfont \makebox(124.0, 82.6)[l]{(d)\strut}}
  \put(1681.05,1208.68){\fontsize{11.76}{14.11}\selectfont \makebox(206.6, 82.6)[l]{$e_1$\strut}}
  \put(1677.85,904.39){\fontsize{11.76}{14.11}\selectfont \makebox(206.6, 82.6)[l]{$e_2$\strut}}
  \put(1440.82,764.41){\fontsize{11.76}{14.11}\selectfont \makebox(206.6, 82.6)[l]{$G_1$\strut}}
  \put(2045.38,1942.00){\fontsize{11.76}{14.11}\selectfont \makebox(124.0, 82.6)[l]{$z$\strut}}
  \put(467.85,562.33){\fontsize{11.76}{14.11}\selectfont \makebox(206.6, 82.6)[l]{$e_1$\strut}}
  \put(1799.72,352.66){\fontsize{11.76}{14.11}\selectfont \makebox(206.6, 82.6)[l]{$e_2$\strut}}
  \put(147.54,270.04){\fontsize{11.76}{14.11}\selectfont \makebox(206.6, 82.6)[l]{$e_3$\strut}}
  \put(423.01,273.24){\fontsize{11.76}{14.11}\selectfont \makebox(206.6, 82.6)[l]{$e_4$\strut}}
  \put(698.48,356.05){\fontsize{11.76}{14.11}\selectfont \makebox(206.6, 82.6)[l]{$e_2$\strut}}
  \put(1582.54,559.12){\fontsize{11.76}{14.11}\selectfont \makebox(206.6, 82.6)[l]{$e_1$\strut}}
  \put(1271.18,263.33){\fontsize{11.76}{14.11}\selectfont \makebox(206.6, 82.6)[l]{$e_3$\strut}}
  \put(1492.85,273.24){\fontsize{11.76}{14.11}\selectfont \makebox(206.6, 82.6)[l]{$e_4$\strut}}
  \put(1624.78,352.66){\fontsize{11.76}{14.11}\selectfont \makebox(206.6, 82.6)[l]{$e_5$\strut}}
  \put(2045.38,1845.23){\fontsize{11.76}{14.11}\selectfont \makebox(124.0, 82.6)[l]{$y$\strut}}
  \put(798.47,1141.75){\fontsize{11.76}{14.11}\selectfont \makebox(124.0, 82.6)[l]{$y$\strut}}
  \put(619.81,1201.30){\fontsize{11.76}{14.11}\selectfont \makebox(124.0, 82.6)[l]{$x$\strut}}
  \put(791.03,1041.25){\fontsize{11.76}{14.11}\selectfont \makebox(124.0, 82.6)[l]{$z$\strut}}
  \put(2026.77,1156.64){\fontsize{11.76}{14.11}\selectfont \makebox(124.0, 82.6)[l]{$y$\strut}}
  \put(2034.06,1048.69){\fontsize{11.76}{14.11}\selectfont \makebox(124.0, 82.6)[l]{$z$\strut}}
  \put(1814.61,1097.08){\fontsize{11.76}{14.11}\selectfont \makebox(124.0, 82.6)[l]{$x$\strut}}
  \end{picture}%
  \caption{The graphs in Subsubcase 2.2.4 of Theorem \ref{thm1.4}.} \label{fig6}
\end{figure}

If $e_1$, $e_2$ are incident with a common vertex in $G_2$, then $G$
must have the structure as given in Figure \ref{fig6} (a).
Similar to the proof of Subsubcase 2.2.1, we can obtain that there
exists an edge cut $F'$ such that $G-F'$ has exactly two components
$G'_1$ and $G'_2$ satisfying that $c(G'_1)=k-1$ if $d_{G_1}(x_2)=1$
or $c(G'_1)=k-2$ if $d_{G_1}(x_2)=2$ and $G'_2$ is a tree of order
$8$. If $G'_1\not\cong K_{2,3}$, then we are done by Claim 1.
Otherwise, $G$ is the graph as given in Figure \ref{fig6} (e) or
(f). In the former case $F''=\{e_1,e_3,e_4\}$ is a good edge cut
while in the latter case $F''=\{e_1,e_3,e_4,e_5\}$ is a good edge
cut.

If $e_1$, $e_2$ are incident with two different vertices in $G_2$,
then $G$ must have the structure as given in Figure \ref{fig6} (b),
(c) or (d). It is easy to see that $F'=\{xy, yz\}$ is a good edge
cut. The proof is thus complete.

{\bf Subsubcase 2.2.5.} $G_1\cong K_{2,3}$ and $G_2\not\cong S_1,
S_3, S_4, W$.

\begin{figure}[ht]
\centering
   \setlength{\unitlength}{0.05 mm}%
  \begin{picture}(2649.7, 1209.7)(0,0)
  \put(0,0){\includegraphics{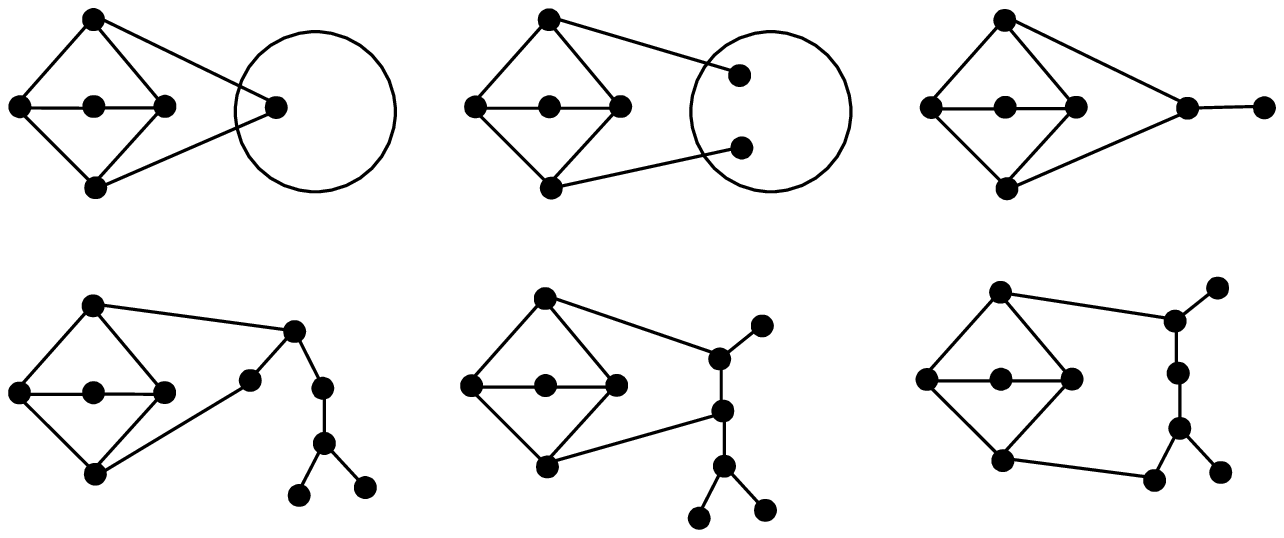}}
  \put(340.97,617.82){\fontsize{11.76}{14.11}\selectfont \makebox(124.0, 82.6)[l]{(a)\strut}}
  \put(264.85,70.02){\fontsize{11.76}{14.11}\selectfont \makebox(124.0, 82.6)[l]{(d)\strut}}
  \put(2174.00,641.20){\fontsize{11.76}{14.11}\selectfont \makebox(124.0, 82.6)[l]{(c)\strut}}
  \put(1277.98,614.66){\fontsize{11.76}{14.11}\selectfont \makebox(124.0, 82.6)[l]{(b)\strut}}
  \put(398.17,1043.94){\fontsize{11.76}{14.11}\selectfont \makebox(206.6, 82.6)[l]{$e_1$\strut}}
  \put(389.11,806.46){\fontsize{11.76}{14.11}\selectfont \makebox(206.6, 82.6)[l]{$e_2$\strut}}
  \put(555.16,675.20){\fontsize{11.76}{14.11}\selectfont \makebox(206.6, 82.6)[l]{$G_2$\strut}}
  \put(44.84,1014.50){\fontsize{11.76}{14.11}\selectfont \makebox(206.6, 82.6)[l]{$e_3$\strut}}
  \put(292.50,1012.00){\fontsize{11.76}{14.11}\selectfont \makebox(206.6, 82.6)[l]{$e_4$\strut}}
  \put(1315.09,1084.62){\fontsize{11.76}{14.11}\selectfont \makebox(206.6, 82.6)[l]{$e_1$\strut}}
  \put(1303.58,758.12){\fontsize{11.76}{14.11}\selectfont \makebox(206.6, 82.6)[l]{$e_2$\strut}}
  \put(1503.61,655.70){\fontsize{11.76}{14.11}\selectfont \makebox(206.6, 82.6)[l]{$G_2$\strut}}
  \put(962.28,1021.34){\fontsize{11.76}{14.11}\selectfont \makebox(206.6, 82.6)[l]{$e_3$\strut}}
  \put(1222.99,1015.98){\fontsize{11.76}{14.11}\selectfont \makebox(206.6, 82.6)[l]{$e_4$\strut}}
  \put(1171.06,37.52){\fontsize{11.76}{14.11}\selectfont \makebox(124.0, 82.6)[l]{(e)\strut}}
  \put(2250.08,1042.36){\fontsize{11.76}{14.11}\selectfont \makebox(206.6, 82.6)[l]{$e_1$\strut}}
  \put(2241.02,804.88){\fontsize{11.76}{14.11}\selectfont \makebox(206.6, 82.6)[l]{$e_2$\strut}}
  \put(1884.97,1013.64){\fontsize{11.76}{14.11}\selectfont \makebox(206.6, 82.6)[l]{$e_3$\strut}}
  \put(2147.62,1007.24){\fontsize{11.76}{14.11}\selectfont \makebox(206.6, 82.6)[l]{$e_4$\strut}}
  \put(55.27,243.13){\fontsize{11.76}{14.11}\selectfont \makebox(206.6, 82.6)[l]{$e_6$\strut}}
  \put(102.90,369.44){\fontsize{11.76}{14.11}\selectfont \makebox(206.6, 82.6)[l]{$e_5$\strut}}
  \put(343.14,529.60){\fontsize{11.76}{14.11}\selectfont \makebox(206.6, 82.6)[l]{$e_1$\strut}}
  \put(359.15,234.91){\fontsize{11.76}{14.11}\selectfont \makebox(206.6, 82.6)[l]{$e_2$\strut}}
  \put(58.06,456.08){\fontsize{11.76}{14.11}\selectfont \makebox(206.6, 82.6)[l]{$e_3$\strut}}
  \put(288.69,427.25){\fontsize{11.76}{14.11}\selectfont \makebox(206.6, 82.6)[l]{$e_4$\strut}}
  \put(970.32,262.54){\fontsize{11.76}{14.11}\selectfont \makebox(206.6, 82.6)[l]{$e_6$\strut}}
  \put(1021.57,384.10){\fontsize{11.76}{14.11}\selectfont \makebox(206.6, 82.6)[l]{$e_5$\strut}}
  \put(1242.58,201.52){\fontsize{11.76}{14.11}\selectfont \makebox(206.6, 82.6)[l]{$e_2$\strut}}
  \put(957.54,474.26){\fontsize{11.76}{14.11}\selectfont \makebox(206.6, 82.6)[l]{$e_3$\strut}}
  \put(1216.96,445.12){\fontsize{11.76}{14.11}\selectfont \makebox(206.6, 82.6)[l]{$e_4$\strut}}
  \put(1283.97,509.34){\fontsize{11.76}{14.11}\selectfont \makebox(206.6, 82.6)[l]{$e_1$\strut}}
  \put(1902.07,278.53){\fontsize{11.76}{14.11}\selectfont \makebox(206.6, 82.6)[l]{$e_6$\strut}}
  \put(1946.91,396.88){\fontsize{11.76}{14.11}\selectfont \makebox(206.6, 82.6)[l]{$e_5$\strut}}
  \put(2116.17,163.38){\fontsize{11.76}{14.11}\selectfont \makebox(206.6, 82.6)[l]{$e_2$\strut}}
  \put(1905.27,473.92){\fontsize{11.76}{14.11}\selectfont \makebox(206.6, 82.6)[l]{$e_3$\strut}}
  \put(2128.27,448.91){\fontsize{11.76}{14.11}\selectfont \makebox(206.6, 82.6)[l]{$e_4$\strut}}
  \put(2209.30,522.13){\fontsize{11.76}{14.11}\selectfont \makebox(206.6, 82.6)[l]{$e_1$\strut}}
  \put(2145.78,40.44){\fontsize{11.76}{14.11}\selectfont \makebox(124.0, 82.6)[l]{(f)\strut}}
  \put(1929.23,823.18){\fontsize{11.76}{14.11}\selectfont \makebox(206.6, 82.6)[l]{$e_6$\strut}}
  \put(1957.68,953.00){\fontsize{11.76}{14.11}\selectfont \makebox(206.6, 82.6)[l]{$e_5$\strut}}
  \end{picture}%
  \caption{The graphs in Subsubcase 2.2.5 of Theorem \ref{thm1.4}.} \label{fig7}
\end{figure}

It is easy to see that $G$ must have the structure as given in
Figure \ref{fig7} (a) or (b). Let $F'=\{e_2,e_3,e_4\}$. Then $G-F'$
has exactly two components $G'_1$ and $G'_2$, where $G'_1$ is a
quadrangle and $G'_2$ is obtained from $G_2$ by adding a pendent
edge. If $G'_2\not\cong S_1, S_3, S_4, W$, then we are done by Claim
1. Otherwise, since $\Delta(G)\leq 3$ and $G_2\not\cong S_1, S_3,
S_4, W$, $G$ must be isomorphic to the graph as given in Figure
\ref{fig3} (c) or Figure \ref{fig7} (c), (d), (e) or (f). In the
first case we are done while in the other cases
$F''=\{e_1,e_4,e_5,e_6\}$ is a good edge cut. The proof is thus
complete.

{\bf Case 3.} $\kappa'(\bar{G})=3$.

Noticing that $\Delta (\bar G)\leq 3$ and $\Delta (G)\leq 3$, we
obtain that $G=\bar G$ is a connected $3$-regular graph.

Let $F=\{e_1, e_2, e_3\}$ be an edge cut of $G$. Then $G-F$ has
exactly two components, say, $G_1$ and $G_2$. Clearly,
$c(G_1)+c(G_2)=k-2\geq 1$.

{\bf Subcase 3.1.} $c(G_1)\geq 1$ and $c(G_2)\geq 1$.

\begin{figure}[ht]
\centering
   \setlength{\unitlength}{0.05 mm}%
  \begin{picture}(1871.2, 771.7)(0,0)
  \put(0,0){\includegraphics{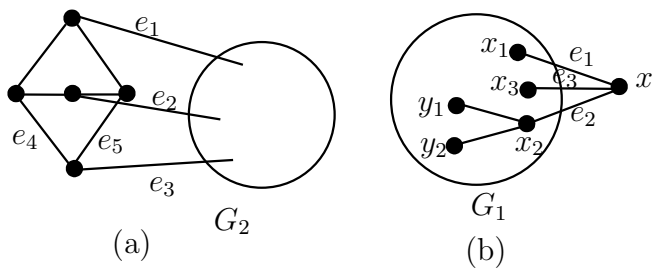}}
  \put(276.61,334.65){\fontsize{11.76}{14.11}\selectfont \makebox(206.6, 82.6)[l]{$e_5$\strut}}
  \put(375.34,646.58){\fontsize{11.76}{14.11}\selectfont \makebox(206.6, 82.6)[l]{$e_1$\strut}}
  \put(413.86,220.70){\fontsize{11.76}{14.11}\selectfont \makebox(206.6, 82.6)[l]{$e_3$\strut}}
  \put(49.09,347.25){\fontsize{11.76}{14.11}\selectfont \makebox(206.6, 82.6)[l]{$e_4$\strut}}
  \put(420.18,451.35){\fontsize{11.76}{14.11}\selectfont \makebox(206.6, 82.6)[l]{$e_2$\strut}}
  \put(585.79,118.97){\fontsize{11.76}{14.11}\selectfont \makebox(206.6, 82.6)[l]{$G_2$\strut}}
  \put(1707.21,497.52){\fontsize{11.76}{14.11}\selectfont \makebox(124.0, 82.6)[l]{$x$\strut}}
  \put(1294.65,590.80){\fontsize{11.76}{14.11}\selectfont \makebox(206.6, 82.6)[l]{$x_1$\strut}}
  \put(1385.92,328.91){\fontsize{11.76}{14.11}\selectfont \makebox(206.6, 82.6)[l]{$x_2$\strut}}
  \put(1530.09,566.16){\fontsize{11.76}{14.11}\selectfont \makebox(206.6, 82.6)[l]{$e_1$\strut}}
  \put(1533.38,406.83){\fontsize{11.76}{14.11}\selectfont \makebox(206.6, 82.6)[l]{$e_2$\strut}}
  \put(1253.69,37.52){\fontsize{11.76}{14.11}\selectfont \makebox(124.0, 82.6)[l]{(b)\strut}}
  \put(1268.85,157.41){\fontsize{11.76}{14.11}\selectfont \makebox(206.6, 82.6)[l]{$G_1$\strut}}
  \put(1320.95,489.65){\fontsize{11.76}{14.11}\selectfont \makebox(206.6, 82.6)[l]{$x_3$\strut}}
  \put(1128.77,435.03){\fontsize{11.76}{14.11}\selectfont \makebox(165.3, 82.6)[l]{$y_1$\strut}}
  \put(1484.25,508.25){\fontsize{11.76}{14.11}\selectfont \makebox(206.6, 82.6)[l]{$e_3$\strut}}
  \put(315.15,55.74){\fontsize{11.76}{14.11}\selectfont \makebox(124.0, 82.6)[l]{(a)\strut}}
  \put(1134.84,323.77){\fontsize{11.76}{14.11}\selectfont \makebox(165.3, 82.6)[l]{$y_2$\strut}}
  \end{picture}%
  \caption{The graphs in Case 3 of Theorem \ref{thm1.4}.} \label{fig8}
\end{figure}

If neither $G_1$ nor $G_2$ is isomorphic to $K_{2,3}$, then we are
done by Claim 1. Otherwise, by symmetry we assume that $G_1\cong
K_{2,3}$. Then $G$ must have the structure as given in Figure
\ref{fig8} (a).
Let $F'=\{e_1,e_2,e_4,e_5\}$. Then it is easy to see that $F'$ is a
good edge cut. The proof is thus complete.

{\bf Subcase 3.2.} One of $G_1$ and $G_2$, say $G_2$ is a tree.

Let $|V(G_2)|=n_2$. Then we have $3n_2=\sum_{v\in
V(G_2)}d_{G}(v)=2(n_2-1)+3=2n_2+1$. Therefore, $n_2=1$, i.e.,
$G_2=S_1$. Let $V(G_2)=\{x\}$, $e_1=xx_1$, $e_2=xx_2$ and
$e_3=xx_3$. Let $N_{G_1}(x_2)=\{y_1, y_2\}$ (see Figure \ref{fig8}
(b)). Let $F'=\{e_1, e_3, x_2y_1, x_2y_2\}$. Then $G-F'$ has exactly
two components $G'_1$ and $G'_2$, where $G'_1$ is a graph obtained
from $G_1$ by deleting a vertex of degree $2$ and $G'_2$ is a tree
of order $2$. Therefore, $c(G'_1)=k-3$. It is easy to check that
$G'_1\not\cong K_{2,3}$. If $G'_1$ is a tree, then we have
$|V(G'_1)|=2$, since $3|V(G'_1)|=\sum_{v\in
V(G'_1)}d_{G}(v)=2(|V(G'_1)|-1)+4=2|V(G'_1)|+2$. Therefore, we are
done by Claim 1. \qed


\begin{thebibliography}{s1}

\bibitem{BM}
J.A. Bondy, U.S.R. Murty, {\it Graph Theory with Applications\/},
Macmillan London and Elsvier, New York (1976).

\bibitem{DS}
J. Day, W. So, Graph energy change due to edge deletion, {\it Lin.
Algebra Appl.\/} 428(2008) 2070--2078.

\bibitem{G1}
I. Gutman, On graphs whose energy exceeds the number of vertices,
{\it Lin. Algebra Appl.\/} 429(2008), 2670--2677.

\bibitem{GLSZ}
I. Gutman, X. Li, Y. Shi, J. Zhang, Hypoenergetic trees, {\it MATCH
Commun. Math. Comput. Chem.\/} 60(2009), 415--426.

\bibitem{GR}
I. Gutman, S. Radenkovi\'c, Hypoenergetic molecular graphs, {\it
Indian J. Chem.\/} 46A (2007), 1733--1736.

\bibitem{LM}
X. Li, H. Ma, Hypoenergetic and strongly hypoenergetic $k$-cyclic
graphs, accepted for publication in {\it MATCH Commun. Math. Comput.
Chem.\/}

\bibitem{MKG}
S. Majstorovi\'{c}, A. Klobu\v{c}ar, I. Gutman, Selected topics from
the theory of graph energy: Hypoenergetic graphs, in: {\it
Applications of Graph Spectra\/}, Math. Inst., Belgrade, 2009,
65--105.

\bibitem{N}
V. Nikiforov, Graphs and matrices with maximal energy, {\it J. Math.
Anal. Appl.\/} 327(2007), 735--738.

\bibitem{N2}
V. Nikiforov, The energy of $C_4$-free graphs of bounded degree,
{\it Lin. Algebra Appl.\/} 428(2008), 2569--2573.

\bibitem{YL}
Z. You, B. Liu, On hypoenergetic unicyclic and bicyclic graphs, {\it
MATCH Commun. Math. Comput. Chem.\/} 61(2009), 479--486.

\bibitem{YLG}
Z. You, B. Liu, I. Gutman, Note on hypoenergetic graphs, {\it MATCH
Commun. Math. Comput. Chem.\/} 62(2009), 491--498.

\end{thebibliography}
\end{document}